\documentclass{amsart}

\usepackage{hyperref}
\usepackage{amssymb, amsmath, amsthm, amsfonts}
\usepackage{verbatim}
\usepackage{bbm}
\usepackage{enumerate}
\usepackage{dsfont}
\usepackage{upgreek}
\usepackage[mathscr]{eucal}
\usepackage{tikz}
\usepackage[normalem]{ulem}

\usepackage[backgroundcolor=black!5,linecolor=black]{todonotes}
\usepackage{comment}

\def\wt#1{\widetilde{#1}}
\def\wh#1{\widehat{#1}}

\newcommand{\N}{\mathbb{N}}
\newcommand{\R}{\mathbb{R}}
\newcommand{\Z}{\mathbb{Z}}

\newcommand{\supp}{\operatorname{supp}}
\newcommand{\dist}{\operatorname{dist}}

\newtheorem{thm}{Theorem}[section]

\newtheorem{prop}[thm]{Proposition}

\numberwithin{equation}{section}

\newcommand\CI{{\mathcal I}}

 \newcommand{\cA}{\mathcal A}

\newcommand{\ep}{{\varepsilon}}

\newcommand{\cP}{\mathcal P}

\newcommand{\Be}{\begin{equation}}
\newcommand{\Ee}{\end{equation}}

\begin{document}

\author[Jeong]{Eunhee Jeong}
\address[Jeong]{Department of Mathematics Education, and  Institute of Pure and Applied Mathematics, Jeonbuk National University, Jeonju 54896, Republic of Korea}
\email{eunhee@jbnu.ac.kr}

\author[Lee]{Sanghyuk Lee}
\address[Lee]{Department of Mathematical Sciences and RIM, Seoul National University, Seoul 08826, Republic of  Korea}
\email{shklee@snu.ac.kr}

\author[Ryu]{Jaehyeon Ryu}
\address[Ryu]{School of Mathematics, Korea Institute for Advanced Study, Seoul
02455, Republic of Korea} 
\email{jhryu@kias.re.kr}
\thanks{This work was supported by NRF (Republic of Korea) grants 2020R1F1A1A01048520, 2022R1A4A1018904, and KIAS Individual Grant MG087001.}

\title[Bochner--Riesz mean for the twisted Laplacian]{
Bochner--Riesz mean for the 
\\
Twisted Laplacian in $\R^2$}

\begin{abstract}
We study the Bochner--Riesz problem for the twisted Laplacian  $\mathcal L$ on $\R^2$.  
For $p\in [1, \infty]\setminus\{2\}$, it has been conjectured that  the Bochner--Riesz means  $S_\lambda^\delta(\mathcal L) f$ of order  $\delta$ converges in $L^p$  for every $f\in L^p$ 
if and only if  $\delta> \max(0,  |(p-2)/p|-1/2)$. We prove  the conjecture by obtaining uniform $L^p$ bounds on $S_\lambda^\delta(\mathcal L)$ up to the sharp summability indices.
\end{abstract}

\maketitle

\section{introduction}
The twisted Laplacian $\mathcal L$ on $\R^{2d}$ is a second order differential operator given  by
\[
\mathcal L = - \sum_{j=1}^d \Big( \big(\frac\partial{\partial x_j}-\frac12iy_j\big)^2 +\big(\frac\partial{\partial{y_j}}+\frac12 ix_j\big)^2\Big),\quad x,y\in \R^d.
\]
The operator $\mathcal L$  is  self-adjoint  and it has  a  discrete spectrum, which is given by the set $2\N_0+d := \big\{2k+d : k\in \N_0\big\}$. Here  $\N_0$ denotes
the set of all natural numbers including $0$.   For $\mu\in 2\N_0+d$, let $\Pi_\mu$ denote the spectral projection operator to the eigenspace with the eigenvalue $\mu$. 
One important property of the projection operators $\Pi_\mu$ is that they allow a spectral decomposition of $L^2$ (\cite{Th93}). That is to say, 
\[ 
f = \sum_{\mu\in 2\N_0+d} \Pi_\mu f, \ \ \ \forall  f\in L^2(\mathbb R^{2d}).
\]

Let $\delta\ge 0$ and $\lambda > 0$. 
By the spectral decomposition,   
the Bochner--Riesz mean $S_\lambda^\delta(\mathcal L)$ for  $\mathcal L$ is defined  by
\[
S_\lambda^\delta(\mathcal L) f = \sum_{\mu\in 2\N_0 + d} \big(1-\frac{\mu}{\lambda}\big)_+^\delta \Pi_\mu f.
\]
The problem known as  the Bochner--Riesz problem is to 
 determine the optimal  summability index $\delta$ for $p\in [1,\infty]$ 
such that $S_\lambda^\delta(\mathcal L) f$ converges to $f$ in $L^p$ for every $f\in L^p$. Of course, this kind of problem was considered  first for 
the Laplacian $-\Delta$ on $\mathbb R^n$, $n\ge 2$, and  the problem  has been extensively studied by numerous authors. It has been conjectured that  the classical Bochner--Riesz mean $S_\lambda^\delta(-\Delta)f$ converges $L^p(\R^n)$ if and only if 
\[
\delta> \delta_\circ(p,n) := \max\Big(0,  n\Big|\frac12-\frac1p\Big|-\frac12 \Big)
\]
for $p\in [1,\infty]\setminus \{2\}$. (When $p=2$, the convergence holds true if and only if $\delta\ge 0$ by  Plancherel's  theorem.)
The conjecture was verified in two dimensions by Carleson and Sj\"olin \cite{CS72}. However, in higher  dimensions, it still remains open and  partial results are known. For the readers who are interested in recent progress on the conjecture, we refer to \cite{Le04, BG11, GHI19, GOWWZ21, Wu23} and references therein.

After the brief digression, we turn back to the Bochner--Riesz problem for $\mathcal L$.
By the uniform boundedness principle, the problem is equivalent to characterizing $\delta=\delta(p)$ such that the estimate
\begin{equation}\label{uniform}
\|S_\lambda^\delta(\mathcal L) f\|_{L^p(\R^{2d})} \le C \|f\|_{L^p(\R^{2d})}
\end{equation}
holds with a constant $C$ independent of $\lambda$ and $f\in \mathcal S(\R^{2d})$.  
In analogue with the classical Bochner--Riesz problem, it is natural to conjecture that 
\eqref{uniform} holds if and only if   $\delta > \delta_\circ(p,2d)$ when $p\neq 2$. 
The necessity part follows by a transplantation theorem due to Kenig--Stanton--Tomas \cite{KST} and the necessary condition for $L^p$ bound on the classical Bochner--Riesz operator  $S_\lambda^\delta(-\Delta)$.

Concerning the sufficiency part,  it was shown by  Thangavelu \cite{Th93} and  Ratnakumar--Rawat--Thangavelu \cite{RRT97}  that \eqref{uniform} holds if $\delta>\delta_\circ(p,2d)$ on a certain  range of $p$.  The range of $p$ was later extended by Stempak and Zienkiewicz \cite{SZ98} for $\max(p,p')>p_\ast(d):={2(2d+1)}/(2d-1)$. All those previous works rely on  a common strategy due to Fefferman and Stein \cite{F2}, which 
 makes it possible  to derive $L^p$ bound on $S_\lambda^\delta(\mathcal L)$  (up to the sharp  exponent $\delta_\circ(p,2d)$)  from the $L^2$--$L^p$ estimate for $\Pi_\mu$: 
\begin{align}\label{e:spectral}
    \|\Pi_\mu f\|_{L^p(\R^{2d})} \le C \mu^{ d(\frac1p-\frac12)-\frac12} \|f\|_{L^2(\R^{2d})}.
\end{align}
The estimate  \eqref{e:spectral} is optimal in that the exponent on $\mu$ cannot be improved.  However, the same strategy does not work any longer if  $\max(p,p')<p_\ast(d)$.  Koch and Ricci \cite{KR07}, in fact, showed that the estimate \eqref{e:spectral} holds if and only if $p_\ast(d)\le p\le \infty$. (See also \cite{JLR} for $L^p$--$L^q$ bounds on $\Pi_\lambda$.)

Other methodologies than the aforementioned have not been exploited until recently in the context of $L^p$ boundedness of $S_\lambda^\delta(\mathcal L)$.  
The second and third named authors \cite{LR22} studied  the problem in a local setting where  $L^p(\mathbb R^{2d})$ is replaced  by $L^p(K)$ for a compact set $K\subset \mathbb R^{2d}$,  and extended the previously known range  for the local $L^p$ bound (\cite{Th98a})  to $\max(p,p')>{2(3d+1)}/(3d-1)$.
 Even though the results are local in their nature, they are more involved than the global bounds on the classical operator $S_\lambda^\delta(-\Delta)$. 
 The local $L^p$ bounds on $S_\lambda^\delta(\mathcal L)$, in fact,  imply the corresponding global bounds on  $S_\lambda^\delta(-\Delta)$ (see \cite{Th98a, LR22}) by virtue of  the transplantation theorem  (\cite{KST}). 
 Remarkably, in $\mathbb R^2$,   the result in \cite{LR22}  gives  the local $L^p$ bounds on the optimal range of $p,\delta$, that is to say, it verifies the Bochner--Riesz conjecture for $\mathcal L$ in a local setting. However, the conjecture without such a local assumption has remained open. 
  
 The objectivity of this article is to prove the Bochner--Riesz conjecture for $\mathcal L$ in $\mathbb R^2$ by obtaining global $L^p$ boundedness of 
 $S_\lambda^\delta(\mathcal L)$.  For the rest of the article, fixing $d=1$, we denote $\delta_\circ(p) := \delta_\circ(p,2)$.  

\begin{thm}\label{thm:ubd}
Let $d=1$ and $1\le p\le \infty$. If $\delta > \delta_\circ(p)$, then the estimate \eqref{uniform} holds.
\end{thm}
 
For a given operator $T$ we denote the kernel of $T$ by $T(z,z')$. 
To prove Theorem \ref{thm:ubd}, we basically  follow the strategy in \cite{LR22} that  is based on kernel expressions 
of the associated multiplier operators (for example, see \eqref{def} below).  The local results  in \cite{LR22} were obtained by combining  asymptotic expansion of the kernel $S_\lambda^\delta(\mathcal L)(z,z')$ 
and estimates for the oscillatory integral operator satisfying {\it Carleson--Sj\"olin} and  ellipticity conditions (\cite{Le06, GHI19}).   More precisely, it was shown that $K_\lambda(z,z'):=S_\lambda^\delta(\mathcal L)(\lambda^{1/2} z, \lambda^{1/2}  z')$ gives rise to an oscillatory integral operator satisfying those conditions under the assumption that   $|z-z'|<2-c$ for a constant $c>0$. 
However, when $(z,z')$ is near the set 
\[\mathfrak S:=\big\{(z,z')\in \mathbb R^2\times \mathbb R^2:  |z-z'|=2\big\}, \] 
the kernel $K_\lambda$ exhibits a different behavior since the critical points of the phase function $\mathcal P(\cdot, z,z')$ (see $\eqref{d:ps}$) are no longer nondegenerate if $(z,z')\in \mathfrak S$. 

To deal with the matter concerning the degeneracy,  we take an  approach inspired by the authors' recent work \cite{JLR2}. 
We make a dyadic decomposition of the kernel away from the set $\mathfrak S$ such that  
the consequent kernels are  supported in the regions $\{(z,z'):  ||z-z'|-2|\sim 2^{-j}\}$. Then, we further break the kernels along  the angle of 
$(z-z')/|z-z'|$ so that each of the decomposed kernels is localized in a set where $z-z'$ is contained in a $ 2^{-j} \times 2^{-j/2}$ rectangle.
Unexpectedly, it turns out that interactions between those angularly decomposed operators are not significant.  
After an appropriate change of variables, we observe  that the operators given by those kernels are  the oscillatory integral operators satisfying the \textit{Carleson--Sj\"olin condition}. 
We combine this observation  with the classical result due to Carleson--Sj\"olin \cite{CS72} to obtain the sharp estimates.

\medskip

\noindent{\it Organization.} 
In Section \ref{sec2}, we break down the proof of Theorem \ref{thm:ubd} to establishing Proposition \ref{prop_main}, which contains the key $L^4$ estimate. 
The subsequent sections  are devoted to proving Proposition \ref{prop_main}. In Section \ref{sec3}  we further reduces the proof so that we only have to  deal with the oscillatory integral operators with kernels supported near $\mathfrak S$ (Proposition \ref{prop:estp}). In Section \ref{sec4}, we complete the proof  by proving  Proposition \ref{prop:estp} via angular decomposition 
and scaling. 

\medskip

\noindent{\it Notations.}
For given non-negative quantities  $A$ and $B$, by $A\lesssim B$ we means that there exists a constant $C>0$ such that $A\le CB.$ We occasionally write $A\lesssim_\ep B$ to indicate that  the implicit constant depends on $\ep>0$.  We write $A\sim B$ if $A\gtrsim B$ and $A\lesssim B$. For an operator $T$,  $\|T\|_{p\to q}$  denotes  the norm of $T$ from $L^p$ to $L^q$.

\section{Reduction to a key $L^4$ estimate}
\label{sec2}
In this section, we make several steps of reduction for the proof of Theorem \ref{thm:ubd} and single out its core part which is Proposition \ref{prop_main} below. 

To prove Theorem \ref{thm:ubd},  it is sufficient to show the estimate \eqref{uniform}  only for $p = 4$ and $\delta>0$. 
Indeed,  the estimate \eqref{uniform} for  two cases $p = 2$, $\delta\ge 0$ and $p = \infty$, $\delta > \delta(\infty)$ are well known (\cite{Th93}). Interpolation with the desired $L^4$ estimate gives \eqref{uniform}  for $2\le p\le \infty$ and $\delta > \delta_\circ(p)$. The case $1\le p <2$ follows by duality.

\subsection{Dyadic decomposition}

Let $\psi \in C^\infty_c([1/4,1])$  such that 
$\sum_{\ell\in\Z} \psi(2^\ell t)=1$ for $t>0.$ 
For $\delta>0$ and $\ell\ge 1$,  set $\psi_\ell^\delta(t)= (2^{-\ell}t)^\delta\psi(2^{-\ell}t)$ and $\psi^\delta_0(t)=t^\delta_+\sum_{\ell\ge0}\psi(2^\ell t)$
so that  
\[ t^\delta_+ =\textstyle \sum_{1\le 2^\ell\le 4\lambda} 2^{\delta\ell}\psi^\delta_\ell(t)\] 
 if $0<t\le \lambda.$ Since $S^\delta_\lambda(\mathcal L)=\lambda^{-\delta}(\lambda-\mathcal L)^\delta_+$, we have 
\[\textstyle S^\delta_\lambda(\mathcal L)=\lambda^{-\delta}\sum_{1\le 2^\ell \le 4\lambda}\,2^{\delta \ell}\psi^\delta_\ell(\lambda-\mathcal L).\]
Therefore, for the estimate \eqref{uniform} for $p = 4$ and $\delta>0$, it is sufficient to show
\begin{equation}\label{e:global}
\|\psi^\delta_\ell(\lambda-\mathcal L)\|_{4\to 4}\lesssim_\ep (\lambda2^{-\ell})^{\ep}, \quad \forall \ep>0.
\end{equation}

By the Fourier inversion, we note 
\begin{align}\label{i:psint}
    \psi^\delta_\ell(\lambda-\mathcal L) = \frac1{2\pi}\int \wh \psi^\delta_\ell(t) e^{i t(\lambda - \mathcal L)} dt.
\end{align}
The kernel of the propagator $e^{-it\mathcal L}$ is given by 
\Be\label{i:epropa}
    e^{-it\mathcal L}(z,z') = c\,(\sin t)^{-1} e^{i(\mathcal P(t,z,z')-t)},  \quad z,z'\in\R^2,
\Ee
 (see \cite{Th93, JLR}) where $c$ is a complex number and 
\begin{align}\label{d:ps}
    \cP(t,z,z'):=t+\frac{|z-z'|^2\cos t}{4\sin t}+\frac{z_2z_1'-z_1z_2'}{2}.
\end{align}

For $\eta\in C^\infty(\R)$, let $[\eta]^\lambda$ be the operator  whose kernel is given by
\begin{equation}\label{def}
[\eta]^\lambda(z,z') = \int \eta(t) (\sin t)^{-1}e^{i\lambda\mathcal P(t,z,z')}dt.
\end{equation}
From \eqref{i:psint} and \eqref{i:epropa}, note that $\psi^\delta_\ell(\lambda-\mathcal L)(\lambda^{1/2}z,\lambda^{1/2} z') = c(2\pi)^{-1}[\wh \psi^\delta_\ell]^\lambda(z,z')$. By scaling, the estimate \eqref{e:global} is equivalent to 
\begin{equation}\label{scaling}
 \|[\wh \psi^\delta_\ell]^\lambda\|_{4\to 4}\lesssim_\ep\lambda^{-1}(\lambda 2^{-\ell})^{\ep},  \quad \forall \ep>0.
\end{equation}

We reduce the proof of \eqref{scaling} to those of the following two propositions.

\begin{prop}\label{prop_main} Suppose that $\eta\in C^\infty_c((0,\pi))$ satisfies $\|{(\frac d{dt})^m}\eta\|_\infty\le C$ for $0\le m\le 100.$ Then, for $\ep>0$ we have
\[ \| [\eta]^\lambda\|_{4\to 4}\lesssim_\ep \lambda^{-1+\ep}.\]
\end{prop}

\begin{prop}\label{lem_origin} Let $1\le 2^\ell \le 4\lambda$ and $\eta\in C^\infty_c((-2^{-5},2^{-5}))$. Then,  for $n\in\Z$  and $\varepsilon>0$ we have
\begin{equation}\label{eq:lemorigin}
 \| [ \eta\wh \psi_\ell^\delta(\cdot -n\pi) ]^\lambda\|_{4\to 4}\lesssim_{\varepsilon} \lambda^{-1}(\lambda 2^{-\ell})^{\ep} (1+2^\ell|n|)^{-9}.
 \end{equation}
\end{prop}

Proposition \ref{prop_main} is the main new contribution of this work, which we prove in the next section, while   Proposition \ref{lem_origin}  is a consequence of the 
local result in  \cite{LR22} (see Proposition \ref{prop:LR} below). Assuming Proposition \ref{prop_main} and \ref{lem_origin} for the moment, we 
prove \eqref{scaling}.

\begin{proof}[Proof of  \eqref{scaling}]
We choose $\eta_0 \in C^\infty_c((-2^{-5},2^{-5}))$ and $\eta_1\in C^\infty_c((2^{-6},\pi-2^{-6}))$ such that both $\eta_0$ and $\eta_1(\cdot +\pi/2)$ are symmetric with respect to $t=0$ and $\eta_0(t)+\eta_1(t)+\eta_0(t-\pi)=1$ for $t\in[0,\pi].$ 
These functions allow us to decompose 
\begin{align}\label{i:decpsi}
    [\wh\psi^\delta_\ell]^\lambda  =\sum_{n\in\Z}\Big( [\eta_0(\cdot +n\pi)\wh\psi^\delta_\ell ]^\lambda +[ \eta_1(\cdot +n\pi)\wh\psi^\delta_\ell]^\lambda
\Big).
\end{align}
Changing variables $t\to t-n\pi$ gives
\[ [ \eta_\kappa(\cdot +n\pi)\wh\psi^\delta_\ell]^\lambda = c[ \eta_\kappa\wh\psi^\delta_\ell(\cdot -n\pi)]^\lambda,\quad \kappa=0,1 \]
with  $|c|=1$. By Proposition \ref{lem_origin} we have 
\begin{equation}\label{e_origin}
 \| [ \eta_0(\cdot +n\pi)\wh \psi_\ell^\delta ]^\lambda\|_{4\to 4}\lesssim_\ep \lambda^{-1}(\lambda 2^{-\ell})^\ep (1+2^\ell|n|)^{-9}.
 \end{equation}
Concerning $[\eta_1\wh\psi^\delta_\ell(\cdot-n\pi)]^\lambda$, note that $\supp \eta_1\subset(0,\pi)$. Since $|\wh\psi^\delta_\ell(t)|\lesssim 
2^\ell(1+2^\ell|t|)^{-M}$ for any $M$,  we have 
$ |(d/{dt})^m(\eta_1(t)\wh\psi^\delta_\ell(t -n\pi))|\lesssim B:= 2^{-\ell M}(1+2^\ell |n|)^{-M}$
 for  $0\le m\le 100$   and any  $M>0$. Applying Proposition \ref{prop_main} to  $\eta=B^{-1} \eta_1\wh\psi^\delta_\ell(\cdot-n\pi)$, we obtain \[ \|[ \eta_1(\cdot +n\pi)\wh\psi^\delta_\ell]^\lambda\|_{4\to 4}\lesssim_\ep \lambda^{-1+\ep}2^{-\ell M}(1+2^\ell|n|)^{-M} \]
 for any  $\ep>0$.   By  \eqref{i:decpsi} and the triangle inequality, using  this and the estimate \eqref{e_origin}, we get  \eqref{scaling}.
\end{proof}

\subsection{Proof of Proposition \ref{lem_origin}}
We use the local estimates for $[\eta]^\lambda$ with a cut-off function $\eta$ supported near the origin, which were  obtained in \cite{LR22}. 

\begin{prop}[{\cite[Theorem 3.3]{LR22}}]\label{prop:LR}
Let $0<\rho<\pi-2^{-5}$ and  $0<c_0<2,$ 
and let $\eta_\rho\in C^\infty_c([2^{-2}\rho,\rho]\cup[-\rho,-2^{-2}\rho])$ satisfy $| (\tfrac d{dt})^m \eta_\rho|\lesssim \rho^{-m}$ for   $0\le m\le 100.$  
Suppose  $E, F\subset \R^2$  be compacts sets such that $|z-z'|\le 2-c_0$ for all $(z, z')\in E\times F$. Then for  $p>4$, we have
\begin{equation}\label{e:LR}
\|\chi_E [\eta_\rho]^\lambda \chi_F\|_{p\to p}\lesssim
 \lambda^{-1}\rho  \,\max\big\{ 1, (\lambda \rho)^{\delta_\circ(p)}\big\}
\end{equation}
\end{prop}

Although the proposition does not include the case $p=4$, interpolation with an easy $L^2$ estimate yields 
\begin{align}\label{e:LR2}
    \|\chi_E [\eta_\rho]^\lambda \chi_F\|_{4\to 4} \lesssim_\ep
\lambda^{-1}\rho \,  \max\big\{ 1, (\lambda\rho)^{\ep}\big\}
\end{align}
for $\ep>0$.  In fact, we note that $\|[\eta_\rho]^\lambda\|_{2\to 2} = \lambda^{-1}\|\int \eta_\rho(t) e^{it(\lambda - \mathcal L)} dt \|_{2\to 2}$.  Thus, it follows that $
\|[\eta_\rho]^\lambda\|_{2\to 2} \le \lambda^{-1}\|\eta_\rho\|_1\lesssim  \lambda^{-1}\rho. 
$ Therefore, $
\| \chi_E [\eta_\rho]^\lambda \chi_F \|_{2\to 2} \lesssim  \lambda^{-1}\rho.
$
Thus, interpolation with \eqref{e:LR} gives \eqref{e:LR2} (taking $p$ arbitrarily close to $4$ when $\rho>\lambda^{-1}$).

We are now ready to prove Proposition \ref{lem_origin}.
Let $\eta \in C^\infty_c((-2^{-5},2^{-5}))$ and let $\mathfrak Q=\{Q\}$ be a tiling of $\R^2$ such that $Q\in\mathfrak Q$ is a square of side length $1/2$. We say $Q\sim Q'$ if $\dist(Q,Q')=0$, and $Q\nsim Q'$ otherwise. Thus, we have  \[ [\eta\wh\psi_\ell^\delta (\cdot-n\pi)]^\lambda = \CI_1+\CI_2,\]
where
\[
\CI_1=\sum_{Q\sim Q'} \chi_Q[\eta\wh\psi_\ell^\delta (\cdot-n\pi)]^\lambda\chi_{Q'}, \qquad \CI_2=\sum_{Q\nsim Q'} \chi_Q[\eta\wh\psi_\ell^\delta (\cdot-n\pi)]^\lambda\chi_{Q'}. \]
The desired estimate \eqref{eq:lemorigin}  follows if we show 
\begin{align}\label{e_local}
 \|\mathcal I_1\|_{4\to 4}&\lesssim \lambda^{-1}(\lambda 2^{-\ell})^{\ep} (1+2^\ell|n|)^{-9}, \\\label{e_error}
\| \CI_2\|_{4\to 4}&\lesssim \lambda^{-N} (1+2^\ell|n|)^{-9}.
\end{align}

To show \eqref{e_local} and \eqref{e_error}, we make an additional decomposition. 
For $j+3\ge 0$, we define
\[
\varphi_j(t) := \psi(2^j t)\quad \text{and} \quad \wt\varphi_j := \psi(2^j|t|).
\]
Then we write
\[
[ \eta \wh \psi^\delta_\ell(\cdot-n\pi) ]^\lambda =\sum_{j\ge 1}[ \eta \wh \psi^\delta_\ell (\cdot -n\pi)\wt\varphi_j ]^\lambda.
\] 
It is easy to see that, for $m\in \N_0$ and any $M$,
\begin{equation}\label{tem1}
\begin{aligned}
 | (\tfrac{d}{dt})^m &(\eta\wh\psi^\delta_\ell (\cdot-n\pi)\wt\varphi_j)|\lesssim   \tilde B:= 2^\ell(1+2^{\ell -j})^{-M} 2^{\max(\ell, j) m} (1+2^\ell|n|)^{-9}.
 \end{aligned}
 \end{equation}
Note that  $|z-z'|\le \sqrt2$ if $(z,z')\in  Q\times Q'$ and $Q\sim Q'$. Thus,  we may apply the estimate \eqref{e:LR2} to $\tilde B^{-1} (\eta\wh\psi^\delta_\ell (\cdot-n\pi)\wt\varphi_j)$.
Using \eqref{tem1} and considering the cases  $\lambda\le 2^j$, $2^\ell\le 2^j<\lambda$, and $2^j\le 2^\ell$, separately,    we obtain 
 \[
\|\chi_Q[\eta \wh\psi^\delta_\ell(\cdot-n\pi)\wt\varphi_j]^\lambda\chi_{Q'}\|_{4\to 4} 
\lesssim \frac{ 2^\ell\lambda^{-1}}{ (1+2^\ell|n|)^{9}} 
\begin{cases}
 \qquad \ 2^{-j},& \lambda\le 2^j,
 \\
\quad \lambda^{\ep}2^{-j(1+\ep)},& 2^\ell\le 2^j<\lambda,
\\
 \lambda^{\ep}2^{-\ell(1+\ep)}2^{(j-\ell)M},& 2^j\le 2^\ell,
 \end{cases}
 \]
 provided that $Q\sim Q'$.
Summation over $j$ yields 
\[ 
\|\chi_Q[ \eta \wh \psi^\delta_\ell (\cdot -n\pi)]^\lambda\chi_{Q'}\|_{4\to 4}\le C \lambda^{-1}(\lambda 2^{-\ell})^{\ep} (1+2^\ell|n|)^{-9}
\] 
with a constant $C$, independent of $\lambda, \ell, n$, and $Q, Q'$ whenever  $Q\sim Q'$. Hence, this gives \eqref{e_local} because for each $Q\in \mathfrak Q$ there are only eight $Q'\in \mathfrak Q$ such that $Q'\sim Q$. 

Now we consider \eqref{e_error}.  
Recall \eqref{def}. The kernel $[ \eta \wh \psi^\delta_\ell (\cdot-n\pi) \wt\varphi_j ]^\lambda(z,z')$ is expressed as an oscillatory integral  with the phase $\mathcal P(t,z,z')$. 
Note that
\begin{align}\label{derivP}
    \partial_t \mathcal P(t,z,z') =1-\frac{|z-z'|^2}{4\sin^2t},
\end{align}
and  $|z-z'|\ge 1/2$ for $(z,z')\in Q\times Q'$ if  $Q\nsim Q'$. Thus,  we have 
\[|\partial_t\mathcal P(t,z,z')|\sim 2^{2j}|z-z'|^2, \quad (t, z,z')\in \supp(\eta\wt\varphi_j)\times Q\times Q'\]
if $Q\nsim Q'$. Moreover, 
$
 |\partial_t^{m+1} \mathcal P(t,z,z')|\lesssim 2^{j(m+1)}|z-z'|^2 $  for any $m\in \mathbb N_0$.
Hence, combining this together with \eqref{tem1}, via repeated integration by parts  we have 
$|[\eta\wh \psi^\delta_\ell(\cdot-n\pi)\wt\varphi_j]^\lambda (z,z')| $ bounded by 
\begin{equation*}\label{tem2}
\begin{aligned}
 K_j(z,z'):=\frac{C}{ (1+2^\ell|n|)^{9}}
\begin{cases}
 2^{M(j-\ell)} (1+\lambda2^{2j}2^{-\ell}|z-z'| )^{-N}, &\ell >j,
 \\[2pt]
2^{\ell -j}  (1+\lambda2^{j}|z-z'|^2)^{-N},  &\ell \le j,
\end{cases}
\end{aligned}
\end{equation*}
for any  $N$ and $M$ if $(z,z')\in Q\times Q'$ and $Q\nsim Q'$.  Thus, the kernel of the operator $\sum_{Q\nsim Q'}\chi_Q[\eta\wh \psi^\delta_\ell (\cdot-n\pi) \wt\varphi_j]^\lambda\chi_{Q'}$ is bounded by  $K_j(z, z')$. Applying Young's inequality, we get
\begin{equation*}\label{Lp_error}
 \|\sum_{Q\nsim Q'}\chi_Q[\eta\wh \psi^\delta_\ell (\cdot-n\pi) \wt\varphi_j]^\lambda\chi_{Q'}\|_{4\to 4}
 \lesssim    \frac{1}{ (1+2^\ell|n|)^{9}}
\begin{cases}
 2^{M(j-\ell)}  (\lambda 2^{-\ell}2^{2j})^{-N}, &\ell >j,\\
2^{\ell-j}  (\lambda 2^{j})^{-N},  &\ell \le j.
\end{cases}
\end{equation*}
for any $N$ and $M$. 
Here, we use the fact that $|z-z'|\gtrsim 1$ if $(z,z')\in Q\times Q'$ and $Q\nsim Q'.$ Hence, \eqref{e_error} follows  by the triangle inequality and summation over $j$.

\section{Dyadic decomposition away from $\mathfrak S$}
\label{sec3}
In this section we prove the key $L^4$ estimate in Proposition \ref{prop_main}.
Throughout this section, we assume that  $\eta\in C_c^\infty((0,\pi))$. 
The first step of the proof is to dyadically decompose the kernel of $[\eta]^\lambda$ near the set $\mathfrak S$.

Recall $\varphi_j=\psi(2^j\cdot)$. 
For $j\in\Z$ satisfying $0\le j \le j_0 := [\log{\lambda^{2/3}}]$, 
we define
\begin{align*}
\chi_j(z,z')&=\varphi_{j-2}(2- |z-z'| ), \\
\chi^\circ(z,z')&=\textstyle \sum_{j> j_0}\wt\varphi_{j-2}(2- |z-z'|).
\end{align*}
Thus, $
 \chi^\circ(z,z')+\sum_{0\le j\le  j_0} \chi_j(z,z')  = 1
$
if  $|z-z'|\le 2$.  We also set
\[
\chi^{e}(z,z') := 1 -\textstyle \big(\sum_{0\le j\le j_0} \chi_j(z,z')  + \chi^\circ(z,z')\big). 
\]
Consequently, $\chi^\circ+\chi^{e}+\sum_{0\le j\le j_0} \chi_j =1$ on $\R^4$. Thus,  
\begin{align}\label{i:suml}
    [\eta]^\lambda = \textstyle \sum_{0 \le j\le j_0} [\eta]^{\lambda}_{j}  + [\eta]^{\lambda,\circ} + [\eta]^{\lambda,e},
\end{align}
where $[\eta]_j^{\lambda}$, $[\eta]^{\lambda, \circ}$, and $[\eta]^{\lambda, e}$ are the operators whose kernels are given by 
\begin{align*}
   [\eta]^{\lambda}_j(z,z') = [\eta]^\lambda(z,z')\chi_{j}(z,z'),  \qquad     [\eta]^{\lambda, \kappa}(z,z') = [\eta]^\lambda(z,z')\chi^{\kappa}(z,z'), \ \kappa \in \{ \circ, e\}.
\end{align*}
As will be seen later, the operators  $[\eta]^{\lambda, \circ}$ and $[\eta]^{\lambda, e}$ are much easier to handle. 

We first prove 
$\|\sum_{0\le j\le j_0} [\eta]_j^{\lambda}\|_{4\to 4}\lesssim_\ep \lambda^{-1+\ep}$
while  the bounds on  the other operators are to be shown near the end of this section.
Since $j_0\lesssim \log\lambda$, 
it suffices to show, for $0\le j\le j_0$, 
\Be 
\label{e:etajgo1}   \|[\eta]_j^{\lambda}\|_{4\to 4}\lesssim_\ep \lambda^{-1+\ep}. 
\Ee

\subsection{Estimate for $[\eta]_j^{\lambda}$}

When $j<C$ for a constant $C$, the desired estimate  \eqref{e:etajgo1} is easy to show by using Proposition \ref{prop:LR}. 
Indeed,  we decompose 
\[
[\eta]_j^{\lambda}= \sum_{-2\le k} [\eta\varphi_k]_j^{\lambda} + \sum_{-2\le k} [\eta\varphi_k(\pi - \cdot)]_j^{\lambda} + [\eta\varphi_0]_j^{\lambda},
\]
where 
$
\varphi_0 = 1 - \sum_{-2\le k} \big(\varphi_k+ \varphi_k(\pi - \cdot)\big).
$
Since $j<C$, $|z-z'|\le 2-c_0$ for a constant $c_0>0$ if $(z,z')\in \supp \chi_j$. 
 By a standard argument (e.g., Proof of Proposition \ref{lem_origin})
 we have
\[
    \|[\eta\varphi_k ]_j^{\lambda} \|_{4\to 4} \lesssim \sup_{B,B'} \|\chi_B [\eta\varphi_k ]_j^{\lambda}\chi_{B'}\|_{4\to 4}
\]
where the supremum is taken over the balls $B,B'$ of radius $c_0/4$
satisfying 
that $\dist(B,B')\le 2- c_0/2$.
Thus, by  \eqref{e:LR2}  we have $\|[\eta\varphi_k ]_j^{\lambda} \|_{4\to 4} \lesssim \lambda^{-1+\varepsilon}.$ 
Moreover,  $[\eta\varphi_k ]_j^{\lambda}= 0$ if $2^k\ge C_\eta$ for a constant $C_\eta>0$.  Therefore, we obtain 
\begin{equation}\label{e:easy_0}
\big\|\sum_{-2\le k} [\eta \varphi_k ]_j^{\lambda} \big\|_{4\to 4} \lesssim_\ep \lambda^{-1+\varepsilon}.
\end{equation}
The same argument  shows $\| [\eta \varphi_0 ]_j^{\lambda} \|_{4\to 4} \lesssim_\ep \lambda^{-1+\varepsilon}.$
To handle $[\eta\varphi_k(\pi - \cdot)]_j^{\lambda} $, we use a symmetric property. 
Considering $\mathbf Lz:=2^{-1/2}(z_1+z_2, z_1-z_2)$, we observe that 
$\mathcal P(\pi-t, \mathbf L z, \mathbf L z') = \pi-\mathcal P(t,z,z')$. Recalling \eqref{def}  and changing  variables $t\to\pi-t$, we see 
\[
[\eta \varphi_k(\pi - \cdot)]_j^{\lambda}(\mathbf L z,\mathbf L z') = C\, \overline{[\eta(\pi-\cdot)\varphi_k]_j^{\lambda}(z,z')}
\]
for a constant $C$ with $|C|=1$. Thus,  $\|[\eta \varphi_k(\pi - \cdot))]_j^{\lambda}\|_{4\to 4} = \|[\eta(\pi-\cdot)\varphi_k]_j^{\lambda}\|_{4\to 4}$. Repeating the previous argument used for \eqref{e:easy_0}, 
we see $\|[\eta \varphi_k(\pi -\, \cdot)]_j^{\lambda}\|_{4\to 4} \lesssim \lambda^{-1+\varepsilon}$. 
Combining all the estimates, we get the bound  \eqref{e:etajgo1} for $j<C$.

Therefore, it is reduced to proving \eqref{e:etajgo1} for $j\ge C$
with a large constant $C$. For the last of this subsection we assume 
$j\ge C$.

\subsubsection*{Further decomposition of the kernel} 
From \eqref{derivP} we  have
\Be
\label{pp}
 | \partial_t \mathcal P(t,z,z')|\sim  |(2-|z-z'|)(2+|z-z'|)- 4\cos^2 t|
 \Ee
 for $(t,z, z')\in \supp \eta\times \supp \chi_j$ since $j\ge C$
with a large constant $C$.  Thus, we are naturally led to decompose dyadically (in $t$) away from 
 $\pi/2$. 

Let $C_0$ be a constant large enough. 
Recalling $\wt\varphi_j= \psi(2^j|\cdot|)$ and $(2-|z-z'|)\sim 2^{-j}$, we decompose 
\[ [\eta]_j^{\lambda} =  [\eta]_{j, 0}^{\lambda} + [\eta]_{j, 1}^{\lambda} \]
where 
\[    [\eta]_{j, 0}^{\lambda}=    \sum_{|2l-j| \le C_0}   [\eta\wt\varphi_l(\pi/2-\cdot)]_j^{\lambda}, \quad [\eta]_{j,1}^{\lambda} = \sum_{|2l-j|>  C_0} [\eta\wt\varphi_l(\pi/2-\cdot)]_j^{\lambda}.\]

\subsubsection*{Estimate for $[\eta]_{j, 1}^{\lambda}$} 
The  estimating $\|[\eta]_{j, 1}^{\lambda}\|_{4\to 4} \lesssim_\varepsilon \lambda^{-1+\varepsilon}$ is easy to obtain.  Indeed, we show this by estimating for the kernel of $ [\eta\wt\varphi_l(\pi/2-\cdot)]_j^{\lambda}$.  
Note that 
\[ [\eta\wt\varphi_l(\pi/2-\cdot)]^{\lambda} (z,z') = \int \mathbf A(t) e^{i\lambda\mathcal P(t,z,z')} dt,  \]
where $\mathbf A(t)= C\eta(t)\wt\varphi_l(\pi/2-t)(\sin t)^{-1}$. 
Since $|2l-j|>C_0$ for a large  $C_0$,  by \eqref{pp}  we get
\[
\big|\lambda \partial_t\mathcal P(t,z,z')\big|\gtrsim  \ \lambda  \max( 2^{-j},  2^{-2l})
\]
for $t\in \supp \mathbf A$ and $(z,z')\in \supp  \chi_j$.  
It is clear that $
    |(d/dt)^n \mathbf A(t)|\lesssim 2^{nl}
    $ for any $ n\in\N_0.$ 
Thus, repeated integration by parts  yields 
\[  |[\eta\wt\varphi_l(\pi/2-\cdot)]_j^{\lambda} (z,z')|\lesssim   b_{l}:=2^{-l}\big(1 + \lambda  2^{-l} \max( 2^{-j},  2^{-2l})\big)^{-N} \] 
for any $N\in\N$. Consequently, we see
$ \sup_z \| [\eta\wt\varphi_l(\pi/2-\,\cdot)]_j^{\lambda} (z,\cdot)\|_1$, $ \sup_{z'} \| [\eta\wt\varphi_l(\pi/2-\,\cdot)]_j^{\lambda} (\cdot, z')\|_1\lesssim 2^{-j} b_l$. 
 Young's inequality and the triangle inequality give
\[   \|[\eta]_{j,1}^{\lambda}\|_{4\to 4} \lesssim  \textstyle \sum_{|2l-j|>  C_0} 2^{-j} b_l\lesssim \lambda^{-1+\varepsilon}. \]

\subsubsection*{Estimate for $[\eta]_{j, 0}^{\lambda}$}
To complete the  proof  of \eqref{e:etajgo1}, it remains to show $\|[\eta\wt\varphi_l(\pi/2-\cdot)]_j^{\lambda}\|_{4\to 4}\lesssim \lambda^{-1+\ep}$ 
 for $1\ll 2^j\le \lambda^{2/3}$ and $2^l\sim 2^{j/2}$.   
Moreover, as before, we note $ \|[\eta \varphi_l(\pi/2-\cdot)]_j^{\lambda}\|_{4\to 4}= \|[\eta\varphi_l(\cdot-\pi/2)]_j^{\lambda}\|_{4\to 4} $  from the symmetric property of the kernel.
Therefore,  the matter is reduced to showing 
\begin{align}\label{e:key}
    \|[\eta\varphi_l(\pi/2-\cdot)]_j^{\lambda}\|_{4\to 4}\lesssim_\ep \lambda^{-1+\ep}
\end{align}
when $1\ll 2^j\le \lambda^{2/3}$ and $2^l\sim 2^{j/2}$. For the purpose, we now consider the stationary point $S_c(z,z')\in (0,\pi/2)$  of the phase function $t\to \mathcal P(t,z,z')$, which is given by 
\Be\label{stpts}
\sin S_c(z,z') = \frac{|z-z'|}{2}.
\Ee
Note that   $\sin (\pi/ 2)- \sin S_c(z,z')\sim 2^{-j}$ if $(z,z')\in \supp \chi_j$. Thus, we have
$S_c(z,z') \in  [ \tfrac\pi 2-c_2 2^{-j/2},   \tfrac\pi 2-c_1 2^{- j/2} ]$
for some $c_1<c_2$ if  $(z,z')\in \supp \chi_j$. 
%

To prove \eqref{e:key}, we make further decomposition of the kernel $[\eta\varphi_l(\pi/2-\cdot)]_j^{\lambda}$
so that $S_c(z,z')$ lies within an interval of length $\ll 2^{-j}$ and the integral for the associated  kernel (for example, \eqref{def}) is also taken over a small interval  of length 
$\ll 2^{-j/2}$.  This can be easily achieved by finite decomposition and discarding some part of the operator which 
has an acceptable $L^4$ bound.

Let $\varepsilon_0>0$ be a sufficiently small constant. 
Recall $\chi_j(z,z')=\psi(2^j(2-|z-z'|))$. 
Breaking $\psi(2^j\cdot)$ into smooth functions supported in finitely overlapping intervals of length $c\varepsilon_0 2^{1-j}$ with a small constant $c>0$,  
we write $\chi_j=\sum  \tilde \chi$ where
\Be 
\label{tchi} \tilde\chi(z,z')= \psi(2^j(2-|z-z'|)) \theta\Big( \frac{a-|z-z'|}{c\varepsilon_0 2^{-j}}\Big)     
\Ee
 for some $a$ satisfying $2-a \in  (2^{-2-j}, 2^{-j}), $ and $\theta\in C_c^\infty((-1,1))$. 
Consequently,  taking $c$ small enough, we have  
\Be
\label{s-pos}
 S_c(z,z')\in  J(t_0, {\varepsilon_0} 2^{-j/2}]:= [t_0-{\varepsilon_0} 2^{-j/2}, t_0+{\varepsilon_0} 2^{-j/2}] 
 \Ee
for some $t_0$  with  $\tfrac \pi 2-t_0 \sim 2^{-j/2}$  if 
$(z, z')  \in \supp \tilde \chi$. Let  $\rho\in C_c^\infty((-2,2))$  such that $\rho=1$ on the interval $[-1,1]$. Set  
\[
\rho_0(t)=\rho\big(2^{j/2}(t-t_0)/{\varepsilon_0}\big). 
\]
Write 
\[  [\eta\varphi_l(\pi/2-\cdot)]^\lambda\tilde\chi=   [\rho_0\eta\varphi_l(\pi/2-\cdot)]^\lambda  \tilde\chi +[(1-\rho_0)\eta\varphi_l(\pi/2-\cdot)]^\lambda\tilde \chi. \]
Here, as before,   $[\eta\varphi_l(\pi/2-\cdot)]^\lambda\tilde\chi$ denotes the operator whose kernel is given by a product of the kernel $[\eta\varphi_l(\pi/2-\cdot)]^\lambda$ and the function  $\tilde \chi$. The other operators are also defined in the same manner.
The operator $[(1-\rho_0)\eta\varphi_l(\pi/2-\cdot)]^\lambda\tilde \chi$ can be easily handled. Indeed, note from \eqref{pp} that $ | \partial_t \mathcal P(t,z,z')|\sim  | \sin S_c(z,z')- \sin t|\sim_{\varepsilon_0} 2^{-j}$ if $(t,z,z')\in \supp (1-\rho_0)\times \supp \chi $. 
Since $2^{2l}\sim 2^j$, by the same argument as before, we have $|[(1-\rho_0)\eta\varphi_l(\pi/2-\cdot)]^\lambda\tilde\chi|\lesssim 2^{- j/2}(\lambda^{-1} 2^{3j/2})^N$. Young's inequality yields $\|[(1-\rho_0)\eta\varphi_l(\pi/2-\cdot)]^\lambda\tilde\chi\|_{4\to 4}\lesssim \lambda^{-1}$.

Therefore, since there are only as many as $O(1/\varepsilon_0)$ $\tilde\chi$, the desired estimate follows if we show  $\|  [\rho_0\eta\varphi_l(\pi/2-\cdot)]^\lambda  \tilde\chi \|_{4\to 4}\lesssim_\varepsilon \lambda^{-1+\varepsilon}$. More generally,  we  prove 
\Be
\label{reduced}
  \|  [\eta]^\lambda \tilde\chi  \|_{4\to 4}\lesssim \lambda^{-1+\varepsilon}
  \Ee
under the following assumption:
\begin{align} 
\label{eta0}
&|(d/dt)^m \eta|\lesssim 2^{mj/2}, \quad \forall m\,;  
\\
\label{eta00}
& \ \  \supp \eta  \subset J(t_0, {\varepsilon_0} 2^{1-j/2}],
\end{align}
for some $t_0$ such that $\tfrac \pi 2-t_0 \sim 2^{-j/2}$. 

\subsubsection*{Asymptotic expansion of the kernel}
We make a change of variables in order that the $t$-derivatives of $\eta$ and $\mathcal P$ are bounded uniformly in $\lambda$ and $j$.
Let us set
\[ \tau (t,z,z') = S_c(z,z') +2^{- j/2}t\]
and 
\begin{align*}
     \tilde{\eta}(t,z,z') = \eta(\tau(t,z,z')), \quad \tilde{\mathcal P}(t,z,z') = 2^{\frac32j}\mathcal P(\tau(t,z,z'),z,z').
\end{align*}

Changing variables  $t\to \tau (t,z,z')$, we have
\[
([\eta]^\lambda \tilde \chi)(z,z') = 2^{-\frac j2} \tilde \chi(z,z')\int \tilde{\eta}(t,z,z')e^{i\lambda 2^{-3j/2}\tilde{\mathcal P}(t,z,z')}dt.
\]
Note that $\supp \tilde \eta(\cdot,z,z')$ is contained in a small interval of length $\lesssim \varepsilon_0$ containing the zero.  We also have 
\[
|\partial_t^m\tilde{\eta}(t,z,z')| \le C_m,  \quad 
|\partial_t^m \big( \tilde{\mathcal P}(t,z,z')- \tilde{\mathcal P}(0,z,z')\big) |\le C_m\]
 for any $m\in \mathbb N_0$ if  $(t,z,z')\in \supp (\tilde\eta \otimes \tilde\chi)$. The former inequality follows from \eqref{eta0}.  
 The latter  inequality for  $m\ge 3$ is clear, and one can show   the inequality for $m=1,2$ using \eqref{pp} and 
\Be
\label{2dl}
\partial_t^2 \mathcal P ( t,z,z') = \frac{|z-z'|^2\cos t }{2\sin^3 t }.
\Ee 
 The case $m=0$ follows from that for $m=1$ via the mean value theorem.  Furthermore,
 since $\partial_t^2\tilde{\mathcal P}(t,z,z')= 2^{j/2}\partial_t^2{\mathcal P}(S_c(z,z')+2^{-j/2}t, z,z')$, from 
 \eqref{2dl} and \eqref{s-pos} we also note that
\[ \partial_t^2\tilde{\mathcal P}(t,z,z')\sim 1,  \quad (t, z,z')\in \supp (\tilde\eta\otimes \tilde\chi).
\]
Since $\partial_t \tilde{\mathcal P}(0,z,z') = 0$, the function $t\to \tilde{\mathcal P}(t,z,z')$ has a nondegenerate critical point at $0$. Taking $\varepsilon_0$ small enough, we apply the stationary phase method 
(\cite[Theorem 7.7.5]{H83}) to obtain the following:  
\begin{align}
\begin{aligned}\label{exp}
([\eta]^\lambda \tilde \chi)(z,z') = \lambda^{-\frac12}2^{\frac j4}\frac{ \tilde \chi(z,z')  \tilde{\eta}(0,z,z')}{ (\partial_t^2\tilde{\mathcal P}(0,z,z')/2\pi)^{1/2}}   e^{i\lambda 2^{-3j/2}\tilde{\mathcal P}(0,z,z')} + E(z,z')
\end{aligned}
\end{align}
where 
 $|E(z,z')|\lesssim \lambda^{-\frac32}2^{\frac 74 j} |\tilde \chi(z,z')|$.
 Note $\|E(\cdot, z')\|_1, \|E(z,\cdot)\|_1\lesssim  \lambda^{-\frac32}2^{\frac 34 j} \lesssim \lambda^{-1}$ since $2^j\lesssim \lambda^{\frac23}$. 
 Thus, Young's inequality shows  
$
\|E\|_{4\to 4}\lesssim \lambda^{-1}.
$

Note that $2^{-3j/2}\tilde{\mathcal P}(0,z,z') = \mathcal P(S_c(z,z'),z,z')$, $\tilde{\eta}(0,z,z')=\eta(S_c(z,z'))$, and $\partial_t^2\tilde{\mathcal P}(0,z,z')=
2^{1+j/2}\cos S_c(z,z')/\sin S_c(z,z')$ using \eqref{2dl} and \eqref{stpts}.  We set 
\begin{align}
\label{ph00}
&\qquad  \quad \Phi(z,z')=  \mathcal P(S_c(z,z'),z,z'),
\\ 
\label{amp} A(z,z')& =  2^{- \frac j4} \tilde\chi(z,z') \eta(S_c(z,z')) \Big(\frac{\sin S_c(z,z')}{\cos S_c(z,z')}\Big)^{1/2}.
\end{align}
For $\mathfrak p\in C^\infty(\mathbb R^4)$ and $\mathfrak a \in C^\infty_c(\mathbb R^4)$, we denote 
 \[  \mathcal T_\lambda [\mathfrak p, \mathfrak a] f(z)= \int e^{i\lambda \mathfrak p (z,z')} \mathfrak a (z,z')   f(z') dz'.\]
Now, by \eqref{exp}  the estimate \eqref{reduced} follows if we show the next proposition.

\begin{prop}\label{prop:estp}
    Let $1\ll 2^j\lesssim \lambda^{2/3}$ and  $\eta$ satisfy \eqref{eta0} and \eqref{eta00}.  Then, $\| \mathcal T_\lambda [\Phi,  A \|_{4\to 4} \lesssim_\ep \lambda^{\ep-1/2} 2^{- j/4}$ for any $\ep>0$. 
\end{prop}

We postpone the proof of Proposition  \ref{prop:estp} until the next section. Before closing this section, we obtain the desired bounds on the operators $[\eta]^{\lambda, \circ}$ and $[\eta]^{\lambda, e}$ (see \eqref{i:suml}).

\subsection{Estimates for  $[\eta]^{\lambda, \circ}$ and $[\eta]^{\lambda, e}$}
In this subsection, we show
\[ \|[\eta]^{\lambda, \circ}\|_{4\to 4}  \lesssim \lambda^{-1}, \quad  \|[\eta]^{\lambda, e}\|_{4\to 4} \lesssim \lambda^{-1}.  \]  
To obtain the above bounds, we use estimates for the kernels.

We first consider $[\eta]^{\lambda, \circ}$. 
Setting $\varphi_\lambda=\sum_{2^{-l}\le C\lambda^{-1/3}} \wt\varphi_l(\pi/2-\cdot)$ for a large positive constant $C$, we decompose 
\begin{align}
\begin{aligned}\label{i:sum2}
    [\eta]^{\lambda, \circ} &= \sum_{2^{-l}> C\lambda^{-1/3}}  [\eta \wt\varphi_l(\pi/2-\cdot)]^{\lambda, \circ} + [\eta\varphi_\lambda]^{\lambda, \circ}.
\end{aligned}
\end{align}
The operator $[\eta\varphi_\lambda]^{\lambda, \circ}$ is easy to handle. Since $|[\eta\varphi_\lambda]^{\lambda, \circ}(z,z')|\lesssim  \lambda^{-1/3}$, it follows that  $\|[\eta\varphi_\lambda]^{\lambda, \circ}(\cdot,z')\|_1,$ $ \|[\eta\varphi_\lambda]^{\lambda, \circ}(z,\cdot)\|_1\lesssim  \lambda^{-1}$.  Consequently, 
 $\|[\eta\varphi_\lambda]^{\lambda, \circ}\|_{4\to 4}\lesssim \lambda^{-1}$.  

As for  $\sum_{2^{-l}> C\lambda^{-1/3}}  [\eta\wt\varphi_l(\pi/2-\cdot)]^{\lambda, \circ} $, we  recall \eqref{def}. Since $2^{-l}> C\lambda^{-1/3}$,  note  from \eqref{pp}  that $ |\partial_t \mathcal P(t,z,z')\big|\gtrsim  2^{-2l}$ if $(t, z,z')\in \supp \wt\varphi_l(\pi/2-\cdot)\times \supp \chi^\circ.$
We also have  $| (d/dt)^n(\eta(t) \wt \varphi_l(\pi/2-t)( \sin t)^{-1})|\lesssim 2^{nl}$ for any $n$. Hence, routine integration by parts gives 
  \[ 
 |[\eta\wt\varphi_l(\pi/2-\cdot)]^{\lambda, \circ} (z,z')|\lesssim 2^{-l}(\lambda^{-1}2^{3l})^N.
 \] 
Thus, we obtain  $ \|[\eta \wt\varphi_l(\pi/2-\cdot)]^{\lambda, \circ}\|_{4\to 4}\lesssim \lambda^{-2/3} 2^{-l}(\lambda^{-1}2^{3l})^N $ by the same argument as before. 
Taking sum over $l$ gives $\|\sum_{2^{-l}> C\lambda^{-1/3}}  [\eta \wt\varphi_l(\pi/2-\cdot)]^{\lambda, \circ} \|_{4\to 4}\lesssim \lambda^{-1}$.

We now turn to  $[\eta]^{\lambda, e}$.  Note $|z-z'|> 2$ for $(z,z')\in \supp \chi^e$. Thus, from \eqref{pp} we have 
$    | \partial_t \mathcal P(t,z,z')|\gtrsim      |z-z'|^2 - 4
    $ for $(t,z,z')\in \supp \eta\times \supp \chi^e.$  Recalling \eqref{def},   by integration by parts  
    we obtain 
\[ 
    |[\eta]^{\lambda, e}(z,z')|\lesssim \big(1+ \lambda(|z-z'|^2-4)\big)^{-N}\chi^e(z,z')
\]
 for any $N\in\N$. Hence, we obtain $\|[\eta]^{\lambda, e}(\cdot,z')\|_1, \|[\eta]^{\lambda, e}(z,\cdot)\|_1\lesssim \lambda^{-1}$. 
   Therefore, we see $\|[\eta]^{\lambda, e}\|_{4\to 4}\lesssim \lambda^{-1}$ by Young's inequality.

\section{$L^4$ bounds  near  the set $\mathfrak S$}
\label{sec4}

In this section,  we prove  Proposition \ref{prop:estp}. 
We begin by decomposing $\mathcal T_\lambda [\Phi, A]$ by breaking 
the amplitude function $A$  along  the angle of $(z-z')/|z-z'|$. 

For each $j$, let  $\Lambda_j\subset \mathbb S^1$ be a collection of $\varepsilon_0 2^{-j/2}$-separated points such that $\mathbb S^1\subset \cup_{\nu\in\Lambda_j} B(\nu, \varepsilon_0 2^{1-j/2})$. 
Let $\{\tilde {\varrho}_j^\nu\}_{\nu\in\Lambda_j}$ be a partition of unity on $\mathbb S^1$ subordinated to $\{B(\nu, \varepsilon_02^{1-j/2})\cap \mathbb S^1\}_{\nu\in \Lambda_j}$. 
Let 
\Be 
\label{an} A^\nu(z,z')= A(z,z') {\varrho}^\nu(z,z'),   \qquad    {\varrho}^\nu(z,z')=\tilde  {\varrho}_j^\nu((z-z')/|z-z'|).
\Ee
Consequently, we  have   $\mathcal T_\lambda [\Phi, A]=\sum_{\nu\in \Lambda_j} \mathcal T_\lambda [\Phi, A^\nu].$
Since $|\Lambda_j|\lesssim 2^{j/2}$,  Proposition \ref{prop:estp} follows once we prove the next. 

\begin{prop}\label{p:l4est0}
    Let $1\ll 2^j \le \lambda^{2/3}$. Then,  for $\nu \in \Lambda_j$ we have
    \Be\label{oj} 
    \| \mathcal T_\lambda [\Phi, A^\nu] \|_{4\to 4}\lesssim  \lambda^{\ep-1/2} 2^{-3j/4}.
    \Ee 
\end{prop}

\subsection{Reduction}
We first make some observations about the operator $\mathcal T_\lambda [\Phi, A]$. 
Note $S_c(z,z')=h(|z-z'|)$ for a function $h$, so  $A(z,z')=a(|z-z'|)$ for a function $a$. 
Thus, the amplitude function $A$ is invariant under simultaneous rotation (i.e., $A(z,z')=A(Rz,Rz')$ for any rotation $R$). 
It is easy to see the phase $\Phi$ is also invariant under  simultaneous rotation. Indeed,  by  \eqref{ph00} and \eqref{d:ps}  we have 
\Be
\label{ph1}  \Phi(z,z')=S_c(z,z')+\cos S_c(z,z')\sin S_c(z,z')+  \mathbf S(z,z'), 
\Ee
where $\mathbf S(z,z')=2^{-1}(z_2 z_1'-z_1 z_2')$. Note $ \mathbf S(z,z')= \mathbf S(Rz,Rz')$ for any notation $R$.  

Therefore, changing variables, it is clear that 
\[\| \mathcal T_\lambda [\Phi, A^\nu] \|_{4\to 4}=\| \mathcal T_\lambda [\Phi(R\cdot, R\cdot), A^\nu(R\cdot, R\cdot)] \|_{4\to 4}.\] 
 As a result, to prove \eqref{oj} we may assume 
\[\nu=e_1.\]

However, due to the term $\mathbf S(z,z')$ in \eqref{ph1},  $\Phi(z,z')$ is not  invariant under simultaneous translation, $(z,z')\to  (z+v,z'+v)$.  
Nevertheless, this does not cause any problem in the perspective of the operator norm. Indeed, 
note that 
\[ \mathbf S(z+v, z'+v')=  \mathbf S(z,z')+ \mathbf S(v, z')+  \mathbf S(z, v')+ \mathbf S(v, v').\]
The second, third, and fourth terms in the phase functions  can be disregarded since they do not have any effect on the operator norm. More generally, we denote 
\[ \Phi_1(z,z')\simeq \Phi_2(z,z')\]
if $\Phi_1(z,z')=\Phi_2(z,z')+ a(z)+ b(z')$ for some functions $a$ and $b$. It is clear that $\| \mathcal T_\lambda [\Phi_1, A^{e_1}] \|_{4\to 4}=\| \mathcal T_\lambda [\Phi_2, A^{e_1}] \|_{4\to 4}$ if 
$\Phi_1(z,z')\simeq \Phi_2(z,z')$.    Using this observation and a standard argument   we can reduce the estimate \eqref{oj}
to a local estimate.

Let $\vartheta\in C_c^\infty((-1,1)^2)$ such that $\sum_{\mathbf k\in \mathbb Z^2} \vartheta(\cdot-\mathbf k)=1$, and set 
$\vartheta_{\mathbf k}(z)= \vartheta( \varepsilon_0^{-1} 2^{j+3} z_1-k_1, \varepsilon_0^{-1} 2^{(j+3)/2}z_2 -k_2)$ for each $\mathbf k =(k_1,k_2)$.
Consequently,  we have 
\[  \mathcal T_\lambda [\Phi,A^{e_1}] f= \textstyle  \sum_{\mathbf k, \mathbf k'} \vartheta_{\mathbf k}  \mathcal T_\lambda [\Phi, A^{e_1}] \vartheta_{\mathbf k'} f.\] 
We also note that 
\Be  
\label{supp} \supp A^{e_1}\subset \big\{ (z_1, z_2):  |z_1-z_1'-a|< c \varepsilon_0 2^{1-j},  \, |z_2-z_2'|< \varepsilon_0^2 2^{-j/2}\,   \big \}   \Ee
with a constant $a$ satisfying $2-a\in (2^{-2-j}, 2^{-j})$.  Thus,  we see that 
$\vartheta_{\mathbf k} \mathcal T_\lambda[\Phi, A^{e_1}] \vartheta_{\mathbf k'} =0$ if  
$|\mathbf k_2-\mathbf k'_2|>3$ and $|\mathbf k_1-  \mathbf k_1'- \varepsilon_0^{-1} 2^{j+3} a|>3 $.
Similarly as in the proof of Proposition \ref{lem_origin}, the estimate \eqref{oj} follows if we show 
\[    \|  \vartheta_{\mathbf k}  \mathcal T_\lambda [\Phi, A^{e_1}] \vartheta_{\mathbf k'}   f\|_4\lesssim_\varepsilon  \lambda^{\ep-1/2} 2^{-3j/4} \|f\|_4\] 
for each $\mathbf k$ and $\mathbf k'$ satisfying $|\mathbf k_2-\mathbf k'_2|\le 2$ and $|\mathbf k_1-  \mathbf k_1'- \varepsilon_0^{-1} 2^{j+3} a|\le 2 $. 
To prove the estimate, thanks to the above discussion, we may use  translation $(z,z')\to (z+v,z'+v)$ for some $v$. 
 Therefore,  we may replace, respectively,  $\vartheta_{\mathbf k}$ and $\vartheta_{\mathbf k'} $ 
 with cutoff functions $\mathfrak a$ and $\mathfrak a'$ such that 
 \Be\label{supp-con}
\begin{aligned}
  & \supp \mathfrak a\subset \{ (z_1, z_2):  |z_1-a| < \varepsilon_0 2^{-j}, \quad   |z_2|< \varepsilon_0 2^{-j/2}\},
   \\
   & \supp \mathfrak a' \subset \{ (z_1, z_2):   |z_1|<\varepsilon_0 2^{-j}, \quad  |z_2|< \varepsilon_0 2^{-j/2}\}, 
   \end{aligned}
   \Ee
 and $\partial^\alpha  \mathfrak a$,  $\partial^\alpha  \mathfrak a'= O(2^{\alpha_2 j}2^{\alpha_2 j/2})$.  Let us set 
 \[ \mathcal A(z,z')= \mathfrak a(z) A^{e_1}(z,z') \mathfrak a'(z').\]
 Therefore, the estimate \eqref{oj} follows from the next.

\begin{prop}\label{p:l4est}
    Let $1\ll 2^j \le \lambda^{2/3}$. Then, $\| \mathcal T_\lambda [\Phi, \mathcal A]\|_{4\to 4} \lesssim_\varepsilon 
    \lambda^{\ep-1/2} 2^{-3j/4}$. 
\end{prop}

Note that $ \mathcal A$ is supported in a product of two rectangle of 
dimension $2^{-j}\times 2^{-j/2}$. We perform change of variables.  Set 
\[  L_j(z,z')= (2^{-j}z_1 +2 , 2^{-j/2} z_2, 2^{-j}z_1',  2^{-j/2} z_2'),\]
and  $\Phi_j=\Phi\circ L_j$ and $\cA_j=\mathcal A\circ L_j$.  Changing variables $(z,z')\to L_j(z,z')$, we have  
\Be 
\label{scaled} \| \mathcal T_\lambda [\Phi, \mathcal A]\|_{4\to 4}=   2^{-3j/2} \|\mathcal T_{ 2^{-3j/2}\lambda} [2^{3j/2}\Phi_j, \cA_j]\|_{4\to 4}.
\Ee
Let $b=2^j(2-a)$ so that  $b\in (2^{-2}, 1)$. From \eqref{supp-con} note that 
\[  \supp \cA_j\subset  U:=  \{(z,z'):  |z_1+b|, |z_1'| < \varepsilon_0, \quad  |z_2|, |z_2'|< \varepsilon_0\}.\]

\subsection{Scaling} 
To obtain estimate for  $\mathcal T_{ 2^{-3j/2}\lambda} [2^{3j/2}\Phi_j, \cA_j]$, we use the known  estimate for the oscillatory integral operator satisfying Carleson--Sj\"olin condition \cite{CS72, H}. For the purpose, we need to take a close look at the scaled  functions
$\cA_j$ and $\Phi_j$. Recalling \eqref{ph00} and \eqref{amp}, we first consider $S_c$ and $|z-z'|$ under  $L_j$. 

Note that  $1-\cos \sigma=g(\sigma^2)$ for an analytic function $g$ with $g(0)=0$ and $g'(0)=1/2$, so    
$g$ has an analytic inverse function near the origin. Consequently, we may write 
$g^{-1}(t)=2t(1+ 2t\mathcal E(2t))$ for an analytic function $\mathcal E$ on a neighborhood of the origin.  
Let us set 
\[ \tilde S(z,z')= \tfrac \pi 2-S_c(z,z') , \qquad \tilde t(z,z')={2- |z-z'|}. \] 
From \eqref{stpts} we have $ 1-\cos \tilde S(z,z')=   \tilde t(z,z')/2 .$  
Recalling that  $|\tilde S(z,z')|\lesssim 2^{-j/2}$, $\tilde t(z,z')\sim 2^{-j}$, and $j>C$ for a large $C$,  from the discussion  above  we have  $ \tilde S^2(z,z')=g^{-1}(\tilde t(z,z')/2)$. Thus, we obtain 
\Be\label{st}
\tilde S(z,z') = \sqrt{\tilde t(z,z')\big(1+ \tilde t(z,z')\mathcal E(\tilde t(z,z'))\big)}.
\Ee
This shows that  $\tilde S$ becomes  singular on the set $\{ (z,z'):\tilde t(z,z')=0\}.$  However, the singularity does not appear  thanks to 
our decomposition. In fact, changing variables $(z,z')\to L_j(z,z')$, we can show the consequent scaled function $2^{j/2} \tilde S\circ L_j$ has bounded derivatives on $U$.  

Indeed,  writing $|z-z'|= (z_1-z_1')(1+{(z_2-z_2')^2}/{(z_1-z_1')^2})^{1/2}$ and using  power series expansion, 
we have  
\[  \tilde t(z,z')=  2-(z_1-z_1')- \frac{(z_2-z_2')^2}{2(z_1-z_1')} \Big(1+O(|z_2-z'_2|^2)\Big).\]
 Thus, it follows that  
\Be 
\label{ttj}
  \tilde t_j(z,z'):= 2^{j} \tilde t(L_j(z,z'))= \mathfrak P(z,z') +O(2^{-j} |z_2-z'_2|^4),   
  \Ee
where 
\[  \mathfrak P(z,z')= z'_1-z_1-   \frac{(z_2-z_2')^2}{2(2+ 2^{-j}(z_1-z_1'))}.\]
In particular,  we note $ \tilde t_j\sim 1$ and $\mathfrak P \sim 1$ on  $U$. 
Combining \eqref{ttj} and \eqref{st} gives 
\Be 
\label{stt}    \tilde S_j(z,z'):= 2^{j/2} \tilde S(L_j(z,z'))= {\mathfrak P}^{\frac12}(z,z')  + \mathcal E(z,z'), 
\Ee 
where  $\mathcal E$ is an analytic function  satisfying 
\Be\label{error} \sup_{(z,z')\in U}  | \partial^\alpha_{z,z'} \mathcal E(z,z')|\lesssim_\alpha 2^{-j}.\Ee

We now claim that 
\Be 
\label{abd}  \sup_{(z,z')\in U} |\partial_{z,z'}^\alpha \cA_j| \le C_\alpha . 
\Ee
For this, it is sufficient to show that the same bound holds for $A^{e_1}\circ L_j$, ${\mathfrak a}_1\circ L_j$, and ${\mathfrak a}_2\circ L_j$.
Those for ${\mathfrak a}_1\circ L_j$ and ${\mathfrak a}_2\circ L_j$ are clear. By \eqref{an} and \eqref{amp}, we need only to show uniform bounds on the derivatives of $\varrho^{e_1}\circ L_j,$
\[\tilde \chi \circ L_j, \quad \mathfrak b_0:=(\eta   \sin^{1/2} )\circ S_c\circ L_j, \quad  \mathfrak b_1:= 2^{-j/4} \cos^{-1/2} \circ\, S_c\circ L_j\] 
  over the set $U$. The bounds on $\varrho^{e_1}\circ L_j$ are easy.  To handle  $\partial^\alpha  \tilde \chi \circ L_j$, we note from \eqref{ttj}  that 
$ \sup_{(z,z')\in U} | \partial^\alpha \tilde t_j(z,z')|\lesssim_\alpha 1$ for any $\alpha$. Thus, from this and  \eqref{tchi} the desired bounds follow. 
For the bounds on $\partial^\alpha \mathfrak b_0$ and $\partial^\alpha \mathfrak b_1$,  by \eqref{stt} we observe that
\Be 
\label{hhh} \sup_{(z,z')\in U} | \partial^\alpha   \tilde S_j (z,z')|\lesssim_\alpha 1.
\Ee
Thus, using \eqref{eta0} and \eqref{eta00}, one can easily see  $ \sup_{(z,z')\in U} | \partial^\alpha \mathfrak b_0(z,z')|\lesssim_\alpha 1$. 
Finally, for $\mathfrak b_1$,  we write  $\mathfrak b_1=  2^{-j/4} \sin^{-1/2} (\tilde S\circ L_j)$  using an elementary trigonometric identity. 
Denote $\varkappa(s)= (s/\sin s)^{1/2}$, which is analytic on $(-\pi, \pi)$.
We write
\[  \mathfrak b_1=2^{-j/4} \sin^{-1/2} (2^{-j/2} \tilde S_j)=   \tilde S_j^{-	1/2} \varkappa (2^{-j/2} \tilde S_j) .\] 
By \eqref{stt} we see $ \tilde S_j\sim 1$ on $U$. Therefore,  $ \sup_{(z,z')\in U} | \partial^\alpha \mathfrak b_1(z,z')|\lesssim_\alpha 1$. 
This proves the claim \eqref{abd}.

We now consider $ p(z,z'):=  S_c+ \cos S_c\sin S_c$. Note that $p(z,z')= \tfrac \pi 2-\tilde S+ 2^{-1} \sin 2\tilde S$. Expanding in power series gives 
\[
 p(z,z') = \tfrac \pi 2-\tfrac23 \tilde S^3 (1 + O(\tilde S^2)). 
\] 
By  \eqref{stt}, we get
\Be 
\label{ph2}
\tfrac \pi 2-  p(L_j(z,z'))= \tfrac23 2^{- \frac{3} 2j } \big({\mathfrak P}^{\frac32}(z,z')+   \mathcal E(z,z') \big ),     
\Ee
where $\mathcal E$ is an analytic  error satisfying \eqref{error}.  Note $\mathbf S\circ L_j(z,z')=2^{-3j/2}\mathbf S(z,z')- 2^{-j/2}z_2'$. By   \eqref{ph1}  we see
\Be  
\label{ph0}  2^{\frac{3} 2j }  \Phi_j(z,z')\simeq  \Phi_j^\ast(z,z')  :=  -\tfrac 23  {\mathfrak P}^{\frac32}(z,z')+    \mathbf S(z,z') + \mathcal E(z,z'),\Ee 
 where $\mathcal E$ is a smooth function satisfying \eqref{error}. 
From \eqref{ph0} it is easy to see that $\sup_{(z,z')\in U} |\partial_{z,z'}^\alpha \Phi_j^\ast| \le C_\alpha $.

\subsection{Carleson--Sj\"olin argument}  To  estimate the right hand side of \eqref{scaled}, we follow the classic argument due to 
Carleson and Sj\"olin \cite{CS72}.  Similarly as before, for $\mathfrak p'\in C^\infty(\mathbb R^3)$ and $\mathfrak a' \in C^\infty_c(\mathbb R^3)$, we denote 
 \[  \mathcal C_\lambda [\mathfrak p', \mathfrak a'] g(z)= \int e^{i\lambda \mathfrak p' (z,s)} \mathfrak a' (z, s)  g(s) ds, \quad (z,s)\in \mathbb R^2\times \mathbb R.\]

Setting  $\Phi_{j}^{\ast, z_1'}(z,s)=\Phi_j^\ast (z,z_1',s)$ and $\cA_{j}^{z_1'}=(z,s)\cA_j (z,z_1',s)$, we observe 
\[\mathcal T_{ 2^{-3j/2}\lambda} [\Phi_j^\ast, \cA_j] f=\int \mathcal C_{2^{-3j/2}\lambda} [\Phi_{j}^{\ast, z_1'}, \cA_{j}^{z_1'}] f(z_1',\cdot) dz'_1.\]
Since $2^j\le \lambda^{2/3}$, thanks to  \eqref{scaled}, 
the desired estimate in Proposition \ref{p:l4est} follows via the Minkowski inequality if we show 
$\|\mathcal C_{\lambda} [\Phi_{j}^{\ast, z_1'}, \cA_{j}^{z_1'}] \|_{4\to4}\lesssim_\ep \lambda^{\ep-1/2}$ for $\lambda\ge 1$.
 For simplicity we make an additional harmless change of variables $z_1\to z_1'-z_1$ so that we can 
replace $\Phi_j^\ast (z,z_1',s)$, $\cA_j (z,z_1',s)$  with 
\[ \Phi_{j, z_1'}^{\ast}(z,s):=\Phi^*_j (z_1'-z_1,z_2,z_1',s), \quad \cA_{j, z_1'}(z,s):=\cA_j (z_1'-z_1,z_2,z_1',s),\] respectively. 
The matter is reduced to showing the uniform bound, for $\lambda\ge 1$, 
\Be
\label{c-s}
 \|\mathcal C_{\lambda} [\Phi_{j, z_1'}^\ast, \cA_{j, z_1'}] \|_{4\to4}\lesssim_\ep \lambda^{\ep-1/2}. 
\Ee

We are now ready to complete the proof of Proposition \ref{p:l4est} by obtaining the estimate \eqref{c-s}. As already mentioned, we use  the well-known 
result regarding the oscillatory integral operator satisfying the Carleson--Sj\"olin condition  \cite{CS72, H}.  The derivatives of $\cA_{j, z_1'}$, $\Phi_{j, z_1'}^\ast$ are, as seen above, uniformly bounded. Thus, for the purpose we only have to show that  $\Phi_{j, z_1'}^{\ast}$ satisfies the Carleson--Sj\"olin condition in a uniform manner. 

Setting  $e(z,s)={(z_2-s)^2}/(2z_1(2- 2^{-j}z_1))$, note 
${\mathfrak P}^{3/2}(z_1'-z_1, z_2,z'_1, s)
=z_1^{3/2} (1-e(z,s)  )^{3/2}$. Expending $(1-e)^{3/2}$ in power series gives
${\mathfrak P}^{3/2}(z_1'-z_1, z_2,z'_1, s)
\simeq -\tfrac32  z_1^{3/2} e(z,s) +  O(|z_2-s|^4)$. Consequently,  we obtain 
\begin{align*} 
{\mathfrak P}^{\frac32}(z_1'-z_1, z_2,z'_1, s)
&\simeq-\tfrac38  {z_1^{1/2}(z_2-s)^2} +  O(|z_2-s|^4) + O(2^{-j}|z_2-s|^2).
\end{align*}
As for $\mathbf S(z, z'_1, s)=2^{-1}(z_2 z_1'-z_1s)$, discarding the harmless term, 
we only need to consider  $-z_1s/2$.  Therefore, recalling \eqref{ph0}, we see
\[ \Phi_{j, z_1'}^{\ast}(z,s)
\simeq   \phi(z,s)
+ \mathcal E(z,z_1', s)+ O(|z_2-s|^4)\]
(here we  abuse the notation $\simeq$) where
\[\phi(z,s)=4^{-1}(2(z_1'-z_1)+2{z_1^{1/2} z_2},\, {z_1^{1/2}} ) \cdot  ( -s, s^2).\]
We note   that  $\phi$ satisfies the Carleson--Sj\"olin condition. Indeed, 
\[   \mathcal M(\phi)(z,s):=\begin{pmatrix}  \nabla_z \partial_s  \phi (z,s)\\   \nabla_z \partial_s^2 \phi(z,s) \end{pmatrix}
 = \frac14
 \begin{pmatrix} 
-1 & 2s 
\\ 
\,\,\,0 & 2
 \end{pmatrix}
 \begin{pmatrix}  -2+ z_1^{-1/2} z_2  &  2 z_1^{1/2}
  \\ 
  z_1^{-1/2}/2  & 0
\end{pmatrix}, 
 \] 
so $\det \mathcal M(\phi)=1/8$.
Recall that 
$\cA_{j,z_1'} (z,s)=0$ unless $|z_1+b|, |z_2|, |s|\le \varepsilon_0$. 
Since $|z_2-s|\le 2\varepsilon_0$ and $j>C$ for a large positive constant $C$, $\Phi_{j, z_1'}^{\ast}$ is a small smooth perturbation of $\phi$,
using stability for $L^p$ bound on the oscillatory integral operator of Carleson--Sj\"olin type,  
 we see that   \eqref{c-s} holds uniformly.

\end{document}